\documentclass[fleqn,a4paper]{article}
\newcommand{\FB}{{\sf FB}}\newcommand{\FIFO}{{\sf FIFO}}
\newcommand{\LAST}{{\sf LAST}}
\newcommand{\pr}{{\sc Proof}}\newcommand{\FBPS}{{\sf FBPS}}
\author{M.F.M. Nuyens\thanks{University of Amsterdam}}
\title{The maximum queue length for heavy-tailed service times
in the M/G/1 FB queue}
\usepackage{amsmath, amsfonts}
\newtheorem{theorem}{Theorem}

\newtheorem{proposition}[theorem]{Proposition}
\newtheorem{corollary}[theorem]{Corollary}
\def\prend{\hfill$\Box$}

\begin{document}
\maketitle
\begin{abstract}\noindent
This paper treats the maximum queue length $M$, in
terms of the number of customers present,
in a busy cycle in the M/G/1 queue.
The distribution of $M$ depends both on the service time distribution and
on the service discipline.
Assume that the service times have a logconvex density and the discipline is 
Foreground Background (\FB).
The \FB\
 service discipline gives service to the customer(s) that have received the least amount of service so far.
It is shown that under these assumptions the tail of $M$ is bounded
by an exponential tail.
This bound is used to calculate the time to overflow
of a buffer, both in stable and unstable queues.
\end{abstract}
\section{Introduction}
In a stochastic process often the extreme values rather than the usual values 
are of great interest. Large values of the queue length may ask for extraordinary measures such as the allocation of auxiliary storage space.
In finite waiting room systems
a natural question is: what is the probability that the queue length will exceed a specific buffer size in a certain time period?  
Unlike the workload process, the queue length process is determined by the  service discipline. 

In this paper we study the maximum queue length $M$, measured in number of customers
present in the system, in a busy period.
The busy period maximum was studied by Cohen \cite{cohen}, who gave an integral representation of the distribution of $M$ in the M/G/1 \FIFO\ queue. For the M/M/1 queue with load $\rho<1$, relation (2.50) in Cohen \cite{cohen} yields a simple expression for the exceedance probabilities:
\begin{equation}\label{ankt} P(M> n)=\rho^n (1-\rho)/(1-\rho^{n+1}).\end{equation}
There is a growing interest in models with
heavy-tailed service times, since statistical data analysis 
has provided convincing evidence of heavy-tailed traffic characteristics
in  high-speed communication networks.
For heavy-tailed service times, the performance of the \FIFO\ discipline decreases and it is natural to consider non-\FIFO\ service disciplines.
 
In this paper the service times have a log-convex\index{log-convex} density $f$, i.e.~$\log f$ is convex.
Both heavy-tailed distributions and
light-tailed distributions may occur. 
Righter \cite{righter} has shown that for queues with log-convex densities
the Foreground Background (\FB) service discipline
minimises the queue length, see Theorem \ref{boot} below. The \FB\ service discipline gives service to the customers that have received the least amount of service so far. 
If there are $n$ such customers, then they are served simultaneously
at rate $1/n$.
When the {\em age}
of a customer is the amount of service he has received, the \FB\ discipline
gives service to the `youngest' customer(s) present in the system. 

Our main result is that in the M/G/1 queue with a log-convex density operating under the \FB\ service discipline 
the tail of the maximum queue length $M$ in a busy period is bounded by an exponential tail:
\begin{equation}\label{mk} P(M> n)\leq \rho^n,\end{equation}  where $\rho$ is the load of the system.

Interestingly, the upper bound does not depend on the precise form of the distribution.
Note furthermore  that (\ref{mk}) is comparable to the precise value of the exceedance probabilities for exponential service
times in (\ref{ankt}). 

By the regenerative structure of the queue length process,
the maximum queue length over a busy period is related to the maximum 
over the time interval $[0,t]$ for $t\to\infty$, see the survey article 
Asmussen \cite{asmussen}.
Using the upper bound (\ref{mk}) we show that in case of service times with heavy-tailed log-convex densities, the time to overflow of a buffer is of another order in the \FB\ queue  than in the \FIFO\ queue. This illustrates the idea that using the \FB\ discipline instead of \FIFO\ may increase the performance of the queue considerably in case of heavy tails.

The class of log-convex densities contains both heavy-tailed distributions like
Pareto distributions (with density $f(x)=(1+x)^{-\alpha}$ for $\alpha>1$) and
Weibull distributions with tail $1-F(x)=\exp(-x^{\beta})$ for $\beta\in(0,1)$ and
light-tailed distributions like the
gamma distribution with density $f(x)=c(\alpha, \beta)x^{\beta}\exp(-\alpha x)$ for $\beta\geq 0$, $\alpha, c(\alpha, \beta)>0$. 
Equivalently, distributions with a log-convex density 
may be characterised by a decreasing likelihood ratio
and have been called DLR distributions in the literature.

Borst {\em et al.~}\cite{borst} and Zwart and Boxma
\cite{zwartboxma} discuss service disciplines in the case of heavy-tailed service times.
For more results on the \FB\ discipline, also known as \FBPS\ or \LAST, see Kleinrock \cite{klein2} and the survey paper Yashkov \cite{yas92}.

The paper is organised as follows. In Section \ref{main} we establish the notation and prove the basic proposition. The inequality
(\ref{mk}) is proved in Section \ref{trm} and a sharper bound depending on 
the Laplace transform of the service-time distribution is established.
In Section \ref{vrr} we use (\ref{mk}) to obtain 
asymptotics for the
maximum queue length $M(t)$ over the interval $[0,t]$ for $t\to\infty$
and apply these to study the overflow time of the buffer in a stable queue.
In Section \ref{vjjf} we describe a coupling that allows us to
apply the results of Section \ref{vrr} to the time to overflow 
in an unstable queue.
\section{Preliminaries}\label{main}
In this section we establish a relation between the maximum queue length
in the M/G/1 \FB\ queue and
 in a related queue. 
This relation is the basis of the theorems in the next section.\\ \\
Consider an M/G/1 \FB\ queue with Poisson arrival rate $\lambda$ and i.i.d.~service times
$B_1, B_2, \ldots$ with distribution function $F$. Let $B$ denote a generic service time and
let the load $\rho=\lambda EB$ satisfy $\rho<1$.
Assume the queue is empty at time 0 and let
$M$ be the maximum queue length in the first busy period. 
In order to determine an upper bound for the tail
probabilities $P(M> n)$, 
we compare the M/G/1 \FB\ queue with a queue 
with the same arrivals and service times, but
with an alternative service discipline, \FB$^*$ say,
defined as follows:
\begin{itemize}
\item
the customer that starts the busy period has the lowest priority,
\index{priority} 
and is only served when there are no other customers present in the queue
\item
the other customers are served according to the \FB\ discipline.
\end{itemize}
Let $M^*$ denote the maximum queue length\index{queue length!maximum} in a busy period in the
M/G/1 \FB$^*$\index{FB@\FB$^{*}$} queue.
The basic idea of this paper is contained in the following proposition
that relates the systems \FB\ and \FB$^*$.
\begin{proposition}\label{honger}
For the discipline \FB$^*$ above
\[ P(M^*> n)=1- Ee^{-\lambda P(M> n-1) B}.\]
\end{proposition}
\pr\
By definition of the discipline \FB$^*$,
the busy period consists of two (disjoint) types of periods:
periods in which only the customer that started the busy period is served and periods during which
only other customers are served. The latter we call {\em sub-busy periods}.
\index{sub-busy periods} Let the random variable $K$ denote the number of sub-busy periods in the busy period.
The sub-busy periods are distributed as M/G/1 \FB\ busy periods with one additional
customer that is not served.
Let the maximum queue lengths in the sub-busy periods be $M_1+1, M_2+1, \ldots, M_K+1$. These are i.i.d.~and have the same distribution as $M+1$.\\
\begin{picture}(350,60)(30,10)
\put(50,27){\line(0,1){6}} 
\put(50,30){\line(1,0){300}}
\put(90,28){\line(0,1){4}}
\put(150,28){\line(0,1){4}}
\put(220,28){\line(0,1){4}}
\put(330,28){\line(0,1){4}}
\put(350,27){\line(0,1){6}}
\put(60,15){$\exp(\lambda)$}
\put(70,35){1}
\put(110,35){$M_1+1$}
\put(185,35){1}
\put(165,15){$\exp(\lambda)$}
\put(250,35){$M_2+1$}
\put(340,35){1}

\end{picture}
\begin{center}Figure 1: 
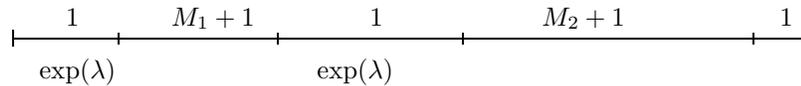
{\em A realisation of the busy period in the M/G/1 \FB$^*$ queue with $K=2$}.\end{center} 
Conditional on the first service time $B_1=x$, 
the number of sub-busy periods $K$ is a Poisson distributed random variable 
with parameter $\lambda x$.
Indeed, $K$ is the number of times the service of the first customer is interrupted.
By the memoryless property of the arrival process, 
the times between two interruptions are independent and exp($\lambda$) distributed.\\  
Given that $K=k$ and $B_1=x$ the maximum queue length $M^*$ is the maximum
of the $k$ independent variables $M_1+1, \ldots, M_k+1$ and hence
\[ P(M^*\leq n\ | \ K=k, B_1=x)=\big(P(M<n)\big)^k.\]
Since $K$ is Poisson distributed for given $x$ we find
\[ P(M^*\leq n, K=k\ | \ B_1=x)=\big(P(M<n)\big)^k\frac{(\lambda x)^k}{k!}e^{-\lambda x}.\]
Summing over $k$ yields a simple expression
\[P(M^*\leq n\ |\ B_1=x)=\sum_{k=0}^{\infty}\big(P(M<n)\big)^k\frac{(\lambda x)^k}{k!}e^{-\lambda x}=
e^{-\lambda xP(M>n-1)}\]
and
\[P(M^*\leq n)= \int_0^{\infty} e^{-\lambda xP(M>n-1)}dF(x)=Ee^{-\lambda P(M>n-1)B}.\]
The result now follows.\prend
\section{Theorems}\label{trm}
In this section we derive exponential 
upper bounds for the tails of the maximum queue length $M$ by combining 
 Proposition \ref{honger} with Theorem \ref{boot} below.\\ \\
The following theorem is Theorem 13.D.8 of Righter \cite{righter}
 in a form adapted to the present setting.
It states that for 
service time distributions with a log-convex density, 
the \FB\ discipline is optimal for minimising the queue length in a strong sense.

As is usual, we restrict attention to service disciplines that cannot look into the future.
This means that each instant the choice of the customer to be served is measurable with respect to the filtration
$\{{\cal F}_t, t\geq 0\}$ where ${\cal F}_t$ is the $\sigma$-algebra
generated by the arrival and departure times and ages of the customers up to time $t$.

\begin{theorem}[Righter]\label{boot}Let $\pi$ be a service discipline and
let  
$N_{\FB}(t)$ and $N_{\pi}(t)$ denote the queue lengths
at time $t$ in G/GI/1 queues with disciplines \FB\ and $\pi$, respectively,
under the same interarrival and service time distribution.
If the service time distribution has a log-convex density, 
then there exist processes $\{\tilde{N}_{\FB}(t), t\geq 0\}$ and
$\{\tilde{N}_{\pi}(t), t\geq 0\}$  such that
\mbox{$\{\tilde{N}_{\FB}(t), t\geq 0\}\stackrel{d}{=}\{N_{\FB}(t),t\geq 0\}$,}
$\{\tilde{N}_{\pi}(t),t\geq 0\}\stackrel{d}{=}\{N_{\pi}(t),t\geq 0\}$ and
\[ P(\tilde{N}_{\FB}(t) \leq \tilde{N}_{\pi}(t), t\geq 0)=1.\]
\end{theorem}
On applying Theorem \ref{boot} to Proposition \ref{honger} it is seen 
that \mbox{$P(M>n)\leq P(M^*>n).$} 
So we find 
\begin{theorem}\label{rsktj}
Let $M$ be the maximum queue length in the busy period 
in the M/G/1 \FB\ queue. If the service time distribution has a log-convex density, 
then the exceedance probabilities $r_n:=P(M>n)$ satisfy
\begin{equation}\label{emktj} r_{n+1}\leq  1- Ee^{-\lambda r_n B},
\qquad n=0, 1,2, \ldots.\end{equation}
\end{theorem}
\begin{theorem}\label{hb} Let $M$ denote the maximum queue length in the busy period in the M/G/1 \FB\ queue with workload $\rho=\lambda EB<1$.
If the service times have a 
log-convex density, then $r_n=P(M>n)$ satisfies 
\begin{equation}\label{nktj} r_n \leq \rho^n,
\qquad n=0,1,\ldots.\end{equation}
In fact,  $r_0=1$, $r_1=1-Ee^{-\lambda B}$ and
\begin{equation}\label{nktj2} r_n\leq \rho\,r_{n-1},
\qquad n\geq 1.\end{equation}
\end{theorem}
\pr\
By Theorem \ref{rsktj} we have $r_{n+1}  \leq \phi (r_n)$,
where $\phi(x)=1-E\exp (-\lambda B x)$.
Furthermore $\phi(0)=0$ and 
\begin{equation}\label{angel} \phi'(x)= E\lambda B \exp(-\lambda x B )\leq \lambda E B=\rho<1.\end{equation}
Hence 
$r_{n+1}\leq \phi(r_n)= \int_0^{r_n} \phi'(u) du\leq
\rho\, r_n$. This is (\ref{nktj2}), which implies
(\ref{nktj}) by induction since $r_0=1$.
The probability $r_1$ may be computed exactly. The first interarrival 
time $U$ is exp($\lambda$)-distributed 
and 
\[ r_1=P(M\geq 2)=P(B_1> U)=1-P(U\geq B_1) =1-\int_0^{\infty} e^{-\lambda x} dF(x)
=1- Ee^{-\lambda B}.\]
This proves the theorem. \hfill$\Box$\\ \\
The approximation
$1-E\exp(-\lambda Bx) \approx \rho x$ in the proof of
Theorem \ref{hb} is quite good for $x$ close to zero, since $\rho$ is the slope of
$\phi(x)=1-E\exp(-\lambda B x)$ in $x=0$.\\ \\
The following
corollary to Theorem \ref{rsktj} gives
even sharper
upper bounds for $P(M> n)$.
\begin{corollary}\label{moei}
Let $q_0=1$ and set $q_{n+1}= 1- Ee^{-\lambda B q_n}$ for $n\geq 1$.
Then
\[ P(M>n)\leq q_n\leq (1-Ee^{-\lambda B})\rho^{n-1},\qquad n\geq 1.\]
\end{corollary}
\pr\ 
By Theorem \ref{rsktj}, induction and the fact that
$1-E\exp(-\lambda B x)$ is increasing in $x$, we have 
 for $n\geq 0$
\[P(M> n) \leq  1- Ee^{-\lambda P(M> n-1) B} 
\leq 1- E\exp(-\lambda q_{n-1} B)=q_n.\]
Inequality (\ref{angel}) and induction yield
$q_{n+1}\leq \rho \, q_n \leq \rho^{n} q_1=(1-Ee^{-\lambda B})\rho^n$.\prend\\ \\
{\bf Example} Numerical calculations show that the bounds $q_n$ in Corollary \ref{moei}
may be significantly better than the bounds in Theorem \ref{hb}.
For service times with the Pareto density
$3(1+x)^{-4}$ and arrival rate \mbox{$\lambda=1.8$}, the load is
$\rho=0.9$ and
 \mbox{$(1-E\exp(-\lambda B))\rho^{99}/q_{100}\approx 7.5$}.\\ \\
If the service times in one queue are stochastically smaller than those
in another queue with the same arrival rate and service discipline, 
the maximum queue length in the first queue
is stochastically smaller than that in the second queue.
Hence the upper bounds in this section also hold 
for service times that are 
stochastically smaller than a service time with a log-convex density.
The next corollary states this idea
for the Pareto distribution.
\begin{corollary}
Consider an M/G/1 \FB\ queue with arrival rate $\lambda$ and generic service time $B$. Let $M$ denote the maximum queue length in a busy period.
Let $\alpha>2$ and $c>0$. If  
\[P(B>x)\leq \frac{1}{(1+cx)^{\alpha}}, 
\qquad x\geq 0,\]  
then  $P(M>n)\leq \theta ^n$ for all $n\geq 0$ where 
$\theta=\lambda /((\alpha-1)(\alpha-2)c^2)$.
\end{corollary}
\section{Asymptotics for the maximum queue length over an interval}\label{vrr}
In this section we present asymptotics for $M(t)$, the maximum queue length
in the interval $(0,t)$, for $t\to\infty$. These are applied to calculate the time to overflow in a finite-buffer system.\\ \\
A stochastic process $X(t)$ is called a {\em regenerative process}\index{regenerative process} if there
exists a (possibly delayed) renewal process\index{renewal process} with epochs $0\leq T_0< T_1
<\cdots$ such that the cycles
\[ \{X(t+T_{i-1})\}_{0\leq t\leq T_i-T_{i-1}}\]
are i.i.d.~for $i=1,2,\ldots$.\\ \\
Consider an M/G/1 \FB\ queue with arrival rate $\lambda$,
i.i.d.~service times with log-convex density $f$, 
generic service time $B$ and workload $\rho<1$. The initial state is arbitrary.
Since the interarrival times are memoryless, the  queue length process is a regenerative process, where a cycle consists of an idle and a busy period. The expected cycle length $\mu$
is given by $\mu=\lambda^{-1}+EB/(1-\rho)=\lambda^{-1}(1-\rho)^{-1}$.
Let $M(t)$ denote the queue length over the time interval $[0,t]$.
Proposition VI.4.7 in Asmussen \cite{asmussenboek} then states that for $t\to\infty$
\begin{equation}\label{maxrelation} \sup_x|P(M(t)\leq x)-P(M\leq x)^{t/\mu}|\to 0.\end{equation}
{\bf Example (Buffer overflow)} We are interested in the time to overflow of a buffer\index{buffer} of size $d\geq 1$.
Denote by 
$t_{d,p}=\inf\{t: P(M(t)> d)\geq p\}$ the lower $p$th quantile for the (random) time to overflow. 
Since $t_{d,p}\to\infty$ as $d\to\infty$, we have by  (\ref{maxrelation}) that 
\[ P(M(t_{d,p})\leq d)=P(M\leq d)^{t_{d,p}/\mu}(1+o(1))\]
as $d\to\infty$.
Theorem \ref{hb} then yields
\begin{equation}\label{boom} t_{d,p}\sim \frac{\mu \log (1-p)}{\log P(M\leq d)} \geq \frac{\mu \log (1-p)}{\log (1-\rho^d)}\end{equation}
where the asymptotic equality holds as $d\to\infty$.
The RHS of (\ref{boom}) may be sharpened by using $q_n$ instead of $\rho^n$, 
where $q_n$ is defined in Corollary \ref{moei}.\\ \\
We compare the time to overflow in the \FB\ and the \FIFO\ queue.
Set $p=\frac{1}{2}$. Then $t_{d,p}$ is the median of the time to overflow. 
Suppose the service times have a Pareto distribution with distribution function $1-(x+1)^{-3}$ and the arrival
rate $\lambda$ is $1.8$, so that the load of this queue is $\rho=0.9$. Time is measured in milliseconds. Let the buffer size $d$ be 1000.
Using the asymptotic inequality (\ref{boom}), we find 
for the \FB\ discipline a value of $t_{d,p}$ larger than $10^{46}$.
This is approximately $10^{35}$ years.

Now consider the same queue with the \FIFO\ discipline.
Occasionally customers with very large service times arrive.
Such customers may cause a buffer overflow.
For large values of $d$, 
the probability that a customer with service time
at least $d/\lambda+3\sigma$, with $\sigma=\sqrt{d/\lambda}$,
 will cause a buffer overflow is larger than 0.99.
The probability that such a customer arrives during the interval
$[0,t]$ is 
\[ p(t)=1-\exp(-\lambda t(1-F\big(d/\lambda+3\sqrt{d/\lambda}\big))).\]
Solving $p(t_1)=1/2$ we find $t_1< 10^8$, approximately one day.

For the M/G/1 \FB\ queue with load $\rho=0.9$ and Weibull service time distribution
 $F(x)=1-\exp(-x^{\beta})$, 
the value of $t_0$ is larger than $10^{48}$ for $\beta=1/4$ and
larger than $10^{46}$ for $\beta=1/2$. For the \FIFO\ queue
with the same service time distribution and load, 
$t_1$ is smaller than $10^3$ for $\beta=1/4$
 and approximately
$10^{20}$ for $\beta=1/2$.
\section{Buffer overflow in unstable queues}\label{vjjf}
In this section we apply the results of section 4 to study the
time to overflow of the buffer in unstable\index{unstable} queues.

When the condition $\rho\leq 1$ is violated, the workload
 asymptotically grows with rate $\rho-1$. Balkema and Verwijmeren \cite{guus} showed that in the {\em overloaded}\index{overloading}  M/G/1 queue under
the \FIFO\ discipline the queue length asymptotically grows linearly with rate $\lambda$, where $\lambda$ is the arrival intensity.
The queue length in the \FB\ queue grows linearly
as well, but with a smaller rate.
By the priority rule of the \FB\ discipline, 
customers with service time less than 
the critical value\index{critical service time} $c^*$, where
$c^*=\inf\{c: \lambda E(B\wedge c)\geq 1\}$, will a.s.~leave the queue,
since they are not hindered by customers with long service times.
Customers with service time larger than $c^*$ have a positive probability of being stuck in the queue forever.
Here we are interested in  the time to overflow of a buffer of size $d$\index{overflow} in an overloaded queue.

Consider a queue with heavy-tailed service times and a small arrival rate.
Assume the stability condition $\rho<1$ is violated. 
This happens because of the mass in the far right tail of the service time distribution. By modifying the right tail, 
we shall obtain an alternative service-time distribution that 
does satisfy the condition $\rho<1$. The maximum queue length in
the alternative queue is comparable to that in the original queue 
for very long periods of time if the arrival process is light.
We shall construct the alternative service times in a way such that they 
have a log-convex density. Then Theorem \ref{hb} and 
(\ref{maxrelation}) are be applied to bound the tail of the maximum queue size.

For the comparison between the two queues we use a coupling argument 
described in Subsection \ref{sap} below. Subsection \ref{vit} illustrates
these ideas with an example.

Extreme value behaviour of $M$ in the M/G/1 \FIFO\ queue was found by Cohen \cite{cohen} under the assumption of the existence of an exponential moment.
Heavy-traffic limits in this case were considered by Iglehart \cite{iglehart}.
Serfozo \cite{serfozo} studied the behaviour of $M(t)$ both in stable and unstable GI/G/1 \FIFO\ queues.

\subsection{The coupling argument}\label{sap}
\begin{theorem}\label{loop}
Let $M_F(t)$ and $M_G(t)$ be the maximum queue lengths over the interval $[0,t]$ 
in two M/G/1 \FB\ queues with service-time distribution functions $F$ and $G$, and
the same arrival rate $\lambda$. Assume $F\wedge p=G\wedge p$ for some
$p\in(0,1)$. Then
\[ M_F(t)\leq_{st}M_G(t)+K_p(t)\]
where $K_p(t)$ is a Poisson process with rate $\lambda (1-p)$.
\end{theorem}
\pr\ Let $N_F(t)$ and $N_G(t)$ be the queue lengths at time $t$. 
We prove that 
\begin{equation}\label{lop} N_F(t)\leq_{st}N_G(t)+K_p(t)\end{equation}
where $K_p(t)$ is a Poisson process with rate $\lambda (1-p)$.
Using the fact that $K_p(t)$ is non-decreasing, the theorem then follows from 
\[M_F(t)=\max_{s\leq t}\{N_F(s)\}\leq_{st} \max_{s\leq t}\{N_G(s)+K_p(s)\}\leq M_G(t)+K_p(t).\]
The coupling is standard:
let $U_1, U_2, \ldots$ be a sequence of i.i.d.~random variables distributed uniformly on $[0,1]$. Define the sequences $B_1, B_2, \ldots$ and $B_1^*, B_2^*, \ldots$ by setting $B_k=F^{-1}(U_k)$ and $B_k^*=G^{-1}(U_k)$ for $k=1, 2,\ldots$.
Then $B_k$ and $B_k^*$ have  distribution functions $F$ and $G$ respectively, for all $k$, and  
\[B_k\wedge p\equiv B_k^*\wedge p,\qquad  k=1,2,\ldots.\]
Now let ${\cal Q}$ and ${\cal Q^*}$ be the M/G/1 \FB\ queues 
with the same arrival process with rate $\lambda$ and
service times introduced above.
Let $N(t)$ and $N^*(t)$ be the queue lengths at time $t$
and let $N_p(t)$ and $N_p^*(t)$ be the part of the queue lengths formed by 
customers younger than $F^{-1}(p)$. 
By the coupling of the service times and the fact that the
\FB\ discipline favours customers younger than $F^{-1}(p)$ over customers older than $F^{-1}(p)$, 
we have $N_p(t)=N_p^*(t)$. 
Let $K_p(t)$ be the Poisson arrival process
of customers with service time larger than $F^{-1}(p)$.
Then $N(t)\leq N_p(t)+K_p(t)$ and $N_p^*(t)\leq N^*(t)$.  Hence
\[ N_F(t)\stackrel{d}{=}N(t)\leq N^*(t)+K_p(t)\stackrel{d}{=} N_G(t)+K_p(t).\]
We conclude the proof by observing that the rate of $K_p(t)$ is $\lambda(1-p)$.\prend\\ \\
Note that the random variables $M_G(t)$ and $K_p(t)$  in Theorem \ref{loop} are dependent.

\subsection{Example}\label{vit}
In this subsection we compute an upper bound for 
 the time to overflow of an M/G/1 \FB\ queue with Pareto
 service density $f(x)=1/(x+1)^2$, so that $\rho=\infty$.
For $a>1$ define
\[ g_a(x)=\begin{cases}
\frac{1}{(x+1)^2} & x\leq a\\
\frac{1}{(a+1)^2} e^{-(x-a)/(a+1)} & x>a.\end{cases}\]
Then $f(x)=g_a(x)$ for $x\leq a$  and it may be shown that $g_a$ is a probability density which is continuous and log-convex.
Now let $M(t)$ and $M_a(t)$ denote the maximum queue lengths over the interval $[0,t]$ of two M/G/1 \FB\ queues with service-time densities $f$ and $g_a$ respectively, and with the same arrival rate $\lambda$.  
Setting $p=1-(a+1)^{-1}$, we have by  Theorem  \ref{loop} 
\begin{equation}\label{bijlk} M(t)\leq_{st} M_a(t)+K_a(t)\end{equation}
where $K_a(t)$ is a Poisson process\index{Poisson!process} with rate 
$\lambda/(1+a)$. Hence for any $x_1,x_2\geq 0$ such that $x_1+x_2=x$ we have
\begin{equation}\label{suitable}P(M(t)>x)\leq P(M_a(t)>x_1)+P(K_a(t)\geq x_2).\end{equation}
Setting $a=10^{40}$ and $\lambda=0.01$, 
we find for the queue with density $g_a$ that
$\rho=\lambda EB=\lambda (\log(a+1)+1)=0.9\ldots$. 
Since for this queue $\rho<1$ and $g_a$ is log-convex, we may apply (\ref{maxrelation}) and (\ref{boom}) to $M_a(t)$.
For a buffer of size $d=1000$ and $p=0.01$ 
we find that $t_{d-1,p}=3\cdot 10^{42}(1+\varepsilon)$, where the error $\varepsilon=\varepsilon(d,p)$ is determined by the asymptotic inequality in (\ref{boom}).
Setting $x_1=d-1$ and $x_2=1$ in (\ref{suitable})
then yields that 
 the probability of a buffer\index{buffer overflow} overflow in $[0,t]$ with $t=10^{40}$ is smaller than $0.02$.

We conclude that for unstable M/G/1 \FB\ queues, the time to overflow may be very large, provided that the arrival rate is low.
\subsection*{Acknowledgements}
I would like to thank A.A.~Balkema for his careful reading of earlier drafts and many useful comments, and W.~Whitt for valuable comments on the presentation of the paper.
Furthermore I am grateful to the referee for his clear remarks and suggestions, which have improved the paper.
\bibliographystyle{plain}\bibliography{bibfile}
\end{document}